\documentclass[10pt,reqno]{amsart} 
\usepackage{amssymb,latexsym,amsmath,amsthm,amsfonts}

\usepackage{color}

\begin{document}

\title[Local Transfer and Reducibility]{Local Transfer and Reducibility of Induced Representations of $p$-adic Groups of Classical Type}
\author[M. Asgari]{Mahdi Asgari} 
\address{Department of Mathematics \\ 
Oklahoma State University \\ 
Stillwater, OK 74078\\
USA} 
\email{asgari@math.okstate.edu} 
\thanks{M.A. was partially supported by Collaboration Grant \# 245422 from the Simons Foundation.}  
	 
\author[J. Cogdell]{James W. Cogdell} 
\address{Department of Mathematics \\ 
The Ohio State University \\ 
Columbus, OH 43210 \\
USA} 
\email{cogdell@math.ohio-sate.edu} 
\thanks{J.W.C. was partially supported by NSF grant DMS--0968505.}  

\author[F. Shahidi]{Freydoon Shahidi} 
\address{Mathematics Department \\  
Purdue University \\ 
West Lafayette, IN 47907 \\ 
USA} 
\email{shahidi@math.purdue.edu} 
\thanks{F.S. was partially supported by NSF grant DMS--1162299.}



\begin{abstract} 
We analyze reducibility points of representations of $p$-adic groups of classical type, 
induced from generic supercuspidal representations of maximal Levi subgroups, both on and 
off the unitary axis.  We are able to give general, uniform results in terms of  local functorial 
transfers of the generic representations of the groups we consider.  The existence of the 
local transfers follows from global generic transfers that were established earlier. 
\end{abstract}

\maketitle

\numberwithin{equation}{section}
\newtheorem{thm}[equation]{Theorem}
\newtheorem{cor}[equation]{Corollary}
\newtheorem{lem}[equation]{Lemma}
\newtheorem{prop}[equation]{Proposition}
\newtheorem{con}[equation]{Conjecture}
\newtheorem{ass}[equation]{Assumption}
\newtheorem{defi}[equation]{Definition}
\newtheorem{exer}[equation]{Exercise}

\theoremstyle{remark}
\newtheorem{rem}[equation]{Remark}
\newtheorem{exam}[equation]{Example}

\newcommand{\Ap}{\A_\infty^+}
\newcommand{\GL}{\operatorname{GL}}
\newcommand{\GO}{\operatorname{GO}}
\newcommand{\GSO}{\operatorname{GSO}}
\newcommand{\GSp}{\operatorname{GSp}}
\newcommand{\GSpin}{\operatorname{GSpin}}
\newcommand{\SL}{\operatorname{SL}}
\newcommand{\SO}{\operatorname{SO}}
\newcommand{\Sp}{\operatorname{Sp}}
\newcommand{\Spin}{\operatorname{Spin}}
\newcommand{\Uef}{\operatorname{U}_{E/F}}

\newcommand{\C}{\mathbb C} 
\newcommand{\Q}{\mathbb Q} 
\newcommand{\Z}{\mathbb Z} 
\newcommand{\R}{\mathbb R} 

\newcommand{\A}{\bf A}
\newcommand{\B}{{\bf B}} 
\newcommand{\G}{{\bf G}} 
\renewcommand{\H}{{\bf H}} 
\newcommand{\M}{{\bf M}} 
\newcommand{\N}{{\bf N}} 
\renewcommand{\P}{{\bf P}} 
\newcommand{\T}{{\bf T}} 
\newcommand{\U}{{\bf U}} 
\newcommand{\X}{{\bf X}} 
\newcommand{\Y}{{\bf Y}} 
 
\newcommand{\sym}{{\rm Sym}}
 
\newcommand{\w}[1]{\ensuremath{\widetilde{#1}}} 
\newcommand{\h}[1]{\ensuremath{\widehat{#1}}}

\section{Introduction}\label{intro} 

In this paper we prove some general, uniform results on reducibility of representation 
induced from irreducible, generic, supercuspidal representations of the Levi subgroups 
of the maximal parabolics of $p$-adic groups.  Some special cases of these results have been 
known for some time.  Our main contribution in this work is to cast these results 
in a general setup in the framework of local Langlands Funcotriality from 
groups of classical type (cf. Section \ref{notation}) to the general linear groups.  
This allows us to prove quite general, uniform 
results using information about poles of local $L$-functions and image of 
local functorial transfers.  Moreover, this agrees with the conjecture on arithmetic $R$-groups, 
defined by Langlands and Arthur \cite[\S 7]{arthur-ast}, through 
Local Langlands Correspondence (cf. \cite{Ban-Goldberg}).

Let $F$ denote a $p$-adic field of characteristic zero.  Consider a maximal Levi 
subgroup of the form $\M=\GL(m) \times \G(n)$ in a larger group $\G(m+n),$ 
a connected linear algebraic group over $F,$  of the same type as $\G.$  
We take the group $\G(n)$ to be any of the split semi-simple groups 
$\SO(2n+1),$ $\Sp(2n),$ $\SO(2n),$ 
the split reductive groups $\GSpin(2n+1),$ $\GSpin(2n),$ 
or the non-split quasi-split groups $\SO_{E/F}(2n),$ $\Uef(2n)$ or $\Uef(2n+1).$ 
Here $E/F$ is a quadratic extension over which our quasi-split group splits.  
These groups all have the property that their connected $L$-group has a classical 
derived group.

Given a connected reductive group $\H$ over $F,$ let $H=\H(F).$  In particular, 
we let $G = \G(F)$ and $M = \M(F),$ where $\G$ and $\M$ are as in the previous 
paragraph. Let $\B  = \T \U$ be a Borel subgroup of $\G$ with $\M \supset \T.$ 
Denote by $\P = \M \N$ the parabolic subgroup of $\G,$ standard via $\N \subset \U$ 
or $\P \supset \B.$  Let $\A_0 \subset \T$ be the maximal split subtorus of $\T$ and 
let $\A \subset \A_0$ be the split component of $\M.$

If $\alpha$ is the unique simple root of $\A_0$ in $\operatorname{Lie}(\N),$ we 
set $\w{\alpha} = \langle \rho, \alpha^\vee \rangle^{-1} \rho,$ where $\rho$ is half the sum 
of the roots of $\A_0$ in $\operatorname{Lie}(\N)$ as in \cite[\S 1.2]{shahidi-book}, $\alpha^\vee$ 
is the coroot of the root $\alpha,$ and $\langle \cdot, \cdot \rangle$ denotes the pairing between 
roots and coroots.  Then, 
$\w{\alpha} \in {\mathfrak a}^*,$ where $\mathfrak a$ is the real Lie algebra of $\A$ and 
$\mathfrak a^*$ is its dual (cf. \cite{shahidi-book}).  Let $s$ be a complex number.  
Then $s \w{\alpha} \in {\mathfrak a}^* \otimes_\R \C.$   Now let $\tau$  be an irreducible 
supercuspidal representation of $M.$  We are interested in understanding reducibility 
of the normalized induced representation 
\begin{equation}
I(s \w{\alpha}, \tau) = 
{\operatorname{Ind}}_{MN}^G 
\left(\tau \otimes q^{\langle s \w{\alpha}, H_M(\cdot) \rangle} \otimes {\bf 1} \right), 
\end{equation}
where $H_M : M \longrightarrow \mathfrak a = 
\operatorname{Hom} \left( X(\M)_F, \R \right)$ 
is defined as in \cite[\S 3.3]{shahidi-book} by 
\begin{equation}
q^{\langle\chi, H_M(m)\rangle} = |\chi(m)|_F, \quad \forall m \in M, 
\end{equation}
where $X(\M)_F$ is the group of the $F$-rational characters of $\M,$  
$\chi \in X(\M)_F$ and $\langle \cdot, \cdot \rangle$ is the pairing between 
$X(\M)_F$ and $\mathfrak a.$

The aim of this paper is to determine the reducibility points for $I(s\w{\alpha}, \tau)$ 
for all $s \in \C,$ whenever $\tau$ is generic, i.e., it has a Whittaker model, in the setting 
of a pair $(\G, \M)$ of classical type, in terms of functoriality as we 
now explain.

By assumption $\M = \GL(m) \times \G(n)$ and in each case ${}^L\G$ embeds 
into $\GL(N,\C) \times W_F$ for a minimal $N,$ with an image with a classical 
derived group (cf. (\ref{HN}) for more detail).  Write $\tau = \sigma \otimes \pi.$  By local transfer, $\pi$ transfers to 
$\Pi = \Pi_1 \boxplus \cdots \boxplus \Pi_d$ on $\GL(N,F)$ (Theorem \ref{local-transfer-sc} here, 
\cite{manuscripta,ckpss1,ckpss2,kk-imrp,kk-aarhus}).

Our main tool is to consider the poles of the intertwining operators (thus zeros of Plancherel measures), as 
proposed by Harish-Chandra \cite{hc-vol4,silberger}, which we determine through poles of certain 
$L$-functions \cite{shahidi1990-annals} (Theorem \ref{const-term} and Corollary \ref{red-cor} here).  
Local transfer allows us to show that these poles exist only when $\sigma$ 
is quasi-self-dual (conjugate-self-dual when $\G(n)$ is unitary) and 
$\sigma$ is of the opposite type to the $L$-group of $\G$ 
(i.e., orthogonal versus symplectic, see Section \ref{loc-red})
or $\sigma$ is among the $\Pi_i$ when it is of the same type as the $L$-group of $\G.$

This provides us with complete information about reducibility points on the unitary 
axis for all groups of classical type. 
The reducibility off the unitary 
axis follows from \cite[Theorem 8.1]{shahidi1990-annals} (Theorem \ref{reduce-off-unitary} here).  
These results are stated as Theorem \ref{reduce-unitary} and summarized for individual groups 
as Propositions \ref{so-odd} -- \ref{unitary}.

The case of induction from other discrete series representation of $M$ must 
also be addressed and it is left for the future.  One should 
also verify the equality of the arithmetic $R$-groups, defined by Langlands and 
Arthur \cite[\S 7]{arthur-ast}, with the analytic $R$-groups, defined by Knapp and Stein \cite{Ban-Goldberg}, 
as conjectured by Langlands and Arthur.

The poles of intertwining operators can also be determined by direct calculations and 
there is a large body of work on this topic, starting with 
\cite{shahidi92-duke,shahidi-invent,goldberg-crelle,GoSh-duke,goldberg-shahidi} 
and ending with some recent work \cite{shahidi-spallone,spallone,Li_w_w-compos,CaiXu-imrn}, 
where the connection to functoriality is fully established in some rather 
general cases ($\SO(2n+1,F)$).

The theory developed in \cite{shahidi1990-annals} (Theorem \ref{reduce-off-unitary} here) applies to any 
quasi-split group and in \cite{lau}, Jing Feng Lau has determined the complete 
reducibility results for exceptional groups $E_6,$ $E_7,$ and  $E_8,$  where 
$M_{\rm der}$ is a product of three $\SL$-groups, using poles of triple product 
$L$-functions which are of Artin type \cite{KimShahidi02-annals}.  The case of exceptional group $G_2$ 
was fully treated in \cite{shahidi1990-annals}

This paper is organized as follows. 
In Section \ref{sec-locrep} we introduce our notation and review the local $L$-functions from the 
Langlands-Shahidi method and their connection to reducibility of induced representations of 
$p$-adic groups of classical type that we consider.  In Section \ref{sc-transfer} we give a proof of the 
generic local transfer of supercuspidal representations of the $p$-adic groups of classical type.  
Our main uniform result on reducibility is given as Theorem \ref{reduce-unitary}. 
The purpose of Section \ref{sec-gps} is then simply to summarize all the information we have, both on and 
off the unitary axis, for each individual group in the hopes that it helps the interested reader see 
what our results give for each individual group and, at the same time, it indicates the scope of these 
results covering all groups of classical type.

{\it Dedication.}  
The first and third named authors would like to dedicate their contributions to this paper to their coauthor, 
Jim Cogdell.  We are very fortunate to have Jim as a friend and collaborator and have very much 
benefited from his kindness and generosity and it is our pleasure to submit this paper to this volume 
in his honor. 

\tableofcontents

\section{Local Representations and $L$-functions} \label{sec-locrep}
 
 \subsection{Notation} \label{notation}
 Let $F$ be a non-archimedean local field of characteristic zero, with $|\cdot|_F$ denoting its 
 $p$-adic absolute value, normalized so that $|\varpi|_F = 1/q,$ where $\varpi$ is a fixed 
 uniformizer of $F$ and $q$ is the the cardinality of the residue field of $F.$  For later use, 
 let us also fix a quadratic extension $E/F.$  Let $\delta_{E/F}$ denote the non-trivial 
 quadratic character of $F^\times$ associated with $E/F$ via Class Field Theory.  We denote 
 the Weil group of $F$ and $E$ by $W_F$ and $W_E,$ respectively. Also, 
 we let $ x \mapsto \bar{x}$ denote the non-trivial element of $\operatorname{Gal}(E/F).$ 
 
 Let $\G$ denote a connected, 
 reductive, quasi-split, linear, algebraic group over $F$. We fix a splitting $(\B, \T, \{X\})$ for 
 $\G$, where $\B$ is a Borel subgroup of $\G$, $\T$ is a maximal  torus in $\B$, and $\{X\}$ is a 
 collection of root vectors, one for each simple root of $\T$ in $\B$. Recall that $\G$ is quasi-split 
 over $F$ if and only if it has an  $F$-splitting, i.e., one preserved under $\operatorname{Gal}(\overline{F}/F)$.

We will assume $\G=\G(n)$ to be one of the following groups:  
the split groups $\SO(2n+1)$, $\Sp(2n)$, $\SO(2n)$, $\GSpin(2n+1)$, $\GSpin(2n),$ 
or the quasi-split non-split groups  
$\Uef(2n)$, $\Uef(2n+1)$, 
$\SO_{E/F}(2n)$, or $\GSpin_{E/F}(2n),$ where $E/F$ is a 
quadratic extension.  
We refer to these groups as {\it groups of classical type}, i.e., those whose connected 
$L$-groups have classical derived groups. 
The groups 
 $\Uef(2n)$ and $\Uef(2n+1)$ are of type ${}^{2}A_{n}$  
and 
 $\SO_{E/F}(2n)$ and $\GSpin_{E/F}(2n)$ are of type ${}^{2}D_n.$  

We write $\B = \T \U$, where $\U$ is the unipotent radical of $\B$.  Unless stated otherwise, 
all the parabolic subgroups we encounter will be assumed to be standard, i.e., they contain $\B.$  

Any standard, maximal, parabolic subgroups $\P$ of $\G$ has a Levi decomposition 
 $\P=\M \N$ with $\M \cong \GL(m) \times \G(n-m),$ if $\G$ is orthogonal or symplectic, 
or  $\M \cong \operatorname{Res}_{E/F} \GL(m) \times \G(n-m),$ 
if $\G$ is unitary.

For later use we define the positive integer $N=N_\G$ to be the dimension of the first fundamental 
representation, or the standard representation, of ${}^{L}\G^{0} = \widehat{\G}(\C)$, 
the connected component of $L$-group of $\G.$  
To be more explicit, for $\G=\G(n)$ as above, we have  
\begin{equation} \label{N}
N = N_\G = \begin{cases}
 2n 
 & 
 \mbox{ if } \G=\Uef(2n), 
\\
 2n+1 
 & 
 \mbox{ if } \G = \Uef(2n+1),  
\\
 2n & \mbox{ if } \G = \SO(2n+1), \GSpin(2n+1), 
 \\
 2n+1 & \mbox{ if } \G = \Sp(2n), 
 \\
2n & \mbox{ if } \G = \SO(2n), \GSpin(2n), 
\\
2n 
& 
\mbox{ if } \G = \SO_{E/F}(2n), \GSpin_{E/F}(2n). 
  	\end{cases}
 \end{equation} 
 In each case the standard representation is a representation of $^{L}G^0$ on $\C^N$ 
 and there is an associated representation of 
 $^{L}G$ on $\C^{N},$ or $\C^{N}\times\C^{N}$ in the unitary group cases,  
 giving rise to a natural $L$-homomorphism 
 \begin{equation}\label{iota}
 \iota : {}^L\G \longrightarrow {}^L\H(N), 
 \end{equation} 
 where 
 \begin{equation} \label{HN}
 \H(N) = 
 \begin{cases}
\GL(N) & \mbox{ if $\G$ is orthogonal or symplectic}, \\
 \operatorname{Res}_{E/F} \GL(N) 
& 
\mbox{ if $\G$ is unitary}. 
 \end{cases}
  \end{equation}
 We refer to \cite[\S1]{cpss-clay} for a detailed description of $\iota$. 
 
 Let $\A_0$ denote the maximal split torus in $\T$ and denote by $\Phi=\Phi(\A_0,\G)$ the restricted 
 roots of $\A_0$ in $\G$ containing positive roots $\Phi^+$. Also, let $\Delta \subset \Phi^+$ denote 
 the set of simple roots. Given a standard maximal parabolic $\P$ there exists a unique 
 $\alpha \in \Delta$ such that $\P=\P_{\theta}$ is determined by the subset 
 $\theta=\Delta \setminus \{\alpha\}$ of $\Delta$. Let $w_0 = w_\G w_\M^{-1}$ denote 
 the longest element of the Weyl group of $\G$ modulo that of $\M$. By abuse of 
 notation, we employ the same symbol for a Weyl group element and its representative in 
 the quotient group. We then have $w_0(\theta) \subset \Delta$ and 
 $w_0(\alpha) < 0$. 

A maximal standard parabolic $\P=\P_\theta$ is called {\it self-associate} if $w_0(\theta) = \theta$. 

\begin{rem} \label{non-self-assoc}
The only non-self-associate case among those we consider above is the following (cf. \cite[\S 4]{kim-fields}): 
\begin{itemize}
\item $D_n$ with $n$ odd and $\theta = \Delta - \{\alpha_n\}$. This corresponds 
to the Levi subgroup $\GL(n)$ in $\SO(2n)$ or $\GL(n)\times\GL(1)$ in $\GSpin(2n)$ with $n$ odd.    
\end{itemize}
\end{rem}

\subsection{The Langlands-Shahidi Local $L$-functions} 

Let $\P=\P_\theta$ be a maximal parabolic in $\G$ as above and let $\rho=\rho_\P$ denote 
half of the sum of positive roots in $\N$. Also, let 
\begin{equation}
\w{\alpha} = \langle \rho, \alpha^\vee \rangle^{-1} \rho. 
\end{equation} 
We have ${}^L\P = {}^L\M {}^L\N$ and 
we let $r$ denote the adjoint action of ${}^L\M$ on the Lie algebra ${}^L\mathfrak n = \operatorname{Lie}({}^L\N).$ 
Let $V_i$ be the subspace of ${}^L\mathfrak n$ spanned by the root vectors $X_{\beta^\vee}$ satisfying 
$\langle \w{\alpha} , \beta^\vee \rangle = i.$ Then we have an irreducible decomposition 
\begin{equation} 
 r = \bigoplus\limits_i r_i, 
\end{equation} 
where $r_i$ denotes the restriction of $r$ to $V_i.$ 

We fix a non-trivial additive character $\psi$ of $F$ throughout.  We can use $\psi$ to define a multiplicative 
character of $\U(F),$ denoted again by $\psi.$  Let $\tau$ be an irreducible $\psi$-generic 
representation of $\M(F)$ and let $s \in \C$. Having fixed $\psi,$ we often simply say $generic$ to mean $\psi$-generic 
in the remainder.  Let 
\begin{equation} 
H_\M : \M(F) \longrightarrow \mathfrak a = \operatorname{Hom}(X(\M)_F,\R) 
\end{equation} 
denote the Harish-Chandra homomorphism defined via 
\begin{equation} 
q^{ \langle \chi, H_\M(m) \rangle } = | \chi(m) |_F , \quad\quad\quad m \in \M(F), \chi \in X(\M)_F. 
\end{equation} 
If $\tau$ is unramified, then it is given by a semisimple conjugacy class $\{ A_\tau\}$ in ${}^L\M$ and 
we set 
\begin{equation} 
L(s,\tau,r_i) = \det\left( I_{V_i} - r_i(\{A_\tau\}) q^{-s} \right)^{-1}. 
\end{equation} 

\subsubsection{Intertwining Operators} 

Let $W = W(\A_0) = N_\G(\A_0)/ \T$ and denote the longest element of $W$ by $w_\ell$. Also, 
let $w_{\ell}^{\M} \in W_\M(\A_0).$ Then $w_0 = w_\ell w_\ell^\M$. Set $\N' = w_0 \N^{-} w_0^{-1}.$  
We define the induced representation 
\begin{equation} \label{ind-s} 
I(s,\tau) = \operatorname{Ind}_{\M(F)\N(F)}^{\G(F)} \left( \tau \otimes q^{\langle s \w{\alpha} , H_\M(\cdot)\rangle} \otimes 1 \right),  
\end{equation}
where $\rho$ denotes half of the sum of positive roots in $\N$ and $\w{\alpha} = \langle \rho, \alpha^\vee \rangle^{-1} \rho.$ 
Here, $\langle \cdot, \cdot \rangle$ denotes the $\mathbb Z$-pairing between characters and cocharacters of $(\G,\T).$  
We also set 
\begin{equation} \label{ind-0}
I(\tau) = I(0,\tau). 
\end{equation}

Define the intertwining operator 
\begin{equation} 
A(s,\tau,w_0) = \int\limits_{\N'(F)} f(w_0^{-1} n' g) dn' : I(s,\tau) \longrightarrow I\left(w_0(s),w_0(\tau)\right). 
\end{equation} 
If $\tau$ is generic, then, via the Langlands-Shahidi method, the $L$-functions $L(s,\tau,r_i)$ are always defined, whether 
$\tau$ is unramified or not, and agree with the definition in the unramified case given above. 

The following two results are well-known (cf. \cite[\S 7]{shahidi1990-annals}).

\begin{thm} \label{const-term}
Assume that $\P$ is a self-associate maximal parabolic and let 
$\tau$ be generic, unitary, supercuspidal. Then 
\[ L(s,\tau,\w{r}_1)^{-1} L(2s,\tau,\w{r}_2)^{-1} A(s,\tau,w_0) \]
is a holomorphic, non-vanishing operator on all of $\C.$
\end{thm}

A consequence of this theorem is the following, which is what we will use later. 

\begin{cor} \label{red-cor}
Suppose that $\tau$ is generic, sucpercuspidal and unitary. 
\begin{itemize}
\item[\bf (a)] If $w_0(\tau) \not\cong \tau$, then $I(\tau)$ is irreducible. (In particular, induction from a 
non-self-associate parabolic is always irreducible.) 
\item[\bf (b)] If $\P$ is self-associate and $w_0(\tau) \cong \tau$, then $I(\tau)$ is irreducible if and only if 
exactly one of $L(s,\tau,\w{r}_1)$ or $L(s,\tau,\w{r}_2)$ has a pole at $s=0$. 
\end{itemize}
(We allow for the second $L$-function not to appear at all.  This does occasionally occur in the 
case of Siegel Levi subgroups, as we will explain later.) 
\end{cor}
 
\section{Generic Local Transfers - Supercuspidal Case} \label{sc-transfer}

One consequence of the generic global functoriality is that we can draw conclusions about transfer of 
local representations, once it is known that the image of the global functorial transfer is an 
isobaric sum of unitary cuspidal representations.  Given that the local transfers are completely 
understood at the archimedean places, we will focus on the non-archimedean local transfers. 

\begin{defi}\label{loc-transfer}
Let $F,$ $\G$ and $N=N_\G$ be as before. 
Let $\pi$ be an irreducible generic representation of $\G(F).$ 
We say an irreducible representation $\Pi$ of $\GL(N,F)$ is a local transfer of $\pi$ if 
\[  L(s,\pi \times \rho) = L(s,\Pi \times \rho) \mbox{ and } \epsilon(s,\pi \times \rho,\psi) = \epsilon(s,\Pi \times \rho, \psi)\]  
or equivalently 
\[  L(s,\pi \times \rho) = L(s,\Pi \times \rho) \mbox{ and } \gamma(s,\pi \times \rho,\psi) = \gamma(s,\Pi \times \rho, \psi) \]  
for all irreducible, unitary, supercuspidal representations $\rho$ of $\GL(m,F)$, $1 \le m \le N-1.$ 
The $L$-, $\epsilon$-, and $\gamma$-factors on the left hand side are those of the Langlands-Shahidi method 
while those on the right hand side are defined via parameters of the Local Langlands Correspondence. 
\end{defi}

We recall that the $\GL \times \GL$ factors on the right hand side are known to equal those defined via 
the Rankin-Selberg or the Langlands-Shahidi methods.

We describe the local transfer for irreducible generic supercuspidal representations in the theorem below.  
This is what we need for our results on reducibility of local representation in Section \ref{loc-red}.
This result is a consequence of the global generic functoriality and was proved in many of 
the cases we cover below along with the global results. We give the details in the proof below.

\begin{thm} \label{local-transfer-sc}

Let $\G=\G(n)$ and $E/F$ be as before. 
Let $\pi$ be an irreducible, generic, supercuspidal, representations of $\G(F).$ 
Then $\pi$ has a unique local transfer $\Pi$ 
to $\GL(N,F)$ if $\G$ is symplectic or orthogonal, 
or to $\GL(N,E)$ if $\G$ is unitary. 
Moreover, it is of the form 
\[ \Pi = \Pi_1 \boxplus \cdots \boxplus \Pi_d = \operatorname{Ind}\left( \Pi_1 \otimes \cdots \otimes \Pi_d \right), \] 
where each $\Pi_i$ is an irreducible, unitary, supercuspidal representation of $\GL(N_i,F)$ 
or $\GL(N_i,E),$ as appropriate, 
and the induction is 
from the standard parabolic subgroup of $\GL(N)$ with Levi component of type $(N_1,\dots,N_d)$ with 
$N_1 + \cdots + N_d = N.$  Furthermore, 

\begin{itemize}

\item if $\G=\SO(2n+1),$ then each $N_i$ is even, each $\Pi_i$ 
satisfies 
$\Pi_i \cong \w{\Pi}_i,$
$L(s, \Pi_i,\wedge^2)$ has 
a pole at $s=0,$ and $\Pi_i \not\cong \Pi_j$ for $i \not= j.$ 

\item if $\G=\SO(2n)$ 
or $\SO_{E/F}(2n),$ 
$n \ge 2,$ or $\G=\Sp(2n),$ then each $\Pi_i$ 
satisfies 
$\Pi_i \cong \w{\Pi}_i,$ 
$L(s, \Pi_i,\sym^2)$ has 
a pole at $s=0,$ and $\Pi_i \not\cong \Pi_j$ for $i \not= j.$ 

\item if $\G=\GSpin(2n+1),$ then each $N_i$ is even, each $\Pi_i$ 
satisfies 
$\Pi_i \cong \w{\Pi}_i \otimes \omega,$ 
where $\omega = \omega_\pi$ is the central character of $\pi,$ 
$L(s, \Pi_i,\wedge^2 \otimes \omega^{-1})$ has 
a pole at $s=0,$ and $\Pi_i \not\cong \Pi_j$ for $i \not= j.$ 

\item if $\G=\GSpin(2n)$ 
or $\GSpin_{E/F}(2n),$  
$n \ge 2,$ then each $\Pi_i$ 
satisfies 
$\Pi_i \cong \w{\Pi}_i \otimes \omega,$ 
where $\omega = \omega_\pi$ is again the central character of $\pi,$ 
$L(s, \Pi_i,\sym^2 \otimes \omega^{-1})$ has 
a pole at $s=0,$ and $\Pi_i \not\cong \Pi_j$ for $i \not= j.$ 

\item if $\G=\Uef(2n+1),$ then each $\Pi_i$ satisfies  
$\overline{\Pi}_i \cong \w{\Pi}_i,$ 
the local Asai $L$-function 
$L(s, \Pi_i, r_A)$ has a pole at $s=0,$ and 
$\Pi_i \not\cong \Pi_j$ 
for $i \not= j.$

\item if $\G=\Uef(2n),$ then each $\Pi_i$ satisfies 
$\overline{\Pi}_i \cong \w{\Pi}_i,$ the local twisted Asai $L$-function 
$L(s, \Pi_i, {r_A \otimes \delta_{E/F}}),$ has a pole at $s=0,$ and 
$\Pi_i \not\cong \Pi_j$ 
for $i \not= j.$  

\end{itemize}  

Here, $\w{\Pi}_i$ denotes the contragredient of $\Pi_i$ and 
$\overline{\Pi}_i$ denotes the 
$\operatorname{Gal}(E/F)$-action on the representation $\Pi,$ 
i.e., $\overline{\Pi}(g) = \Pi(\bar{g}).$ 

\end{thm} 

\begin{proof}
For $\G=\SO(2n+1),$ $\Sp(2n),$ or $\SO(2n),$ this is 
\cite[Theorem 7.3]{ckpss2}.  For $\G=\Uef(2n)$ this is 
\cite[Proposition 8.4]{kk-imrp} and for $\G=\Uef(2n+1)$ 
it is \cite[Proposition 4]{kk-aarhus}.  For $\G=\GSpin(2n+1)$ 
or $\GSpin(2n)$ this is \cite[Theorem 4.7]{hk}.   All of these 
results are based on a local-global argument, using the 
fact that one can embed the generic supercuspidal 
representation $\pi$ as the local component of a global 
generic representation to which one can apply the global 
generic transfer,  possibly several times, to obtain the result.  
Let us give some of the details now to indicate that a similar 
argument works for all the groups we are considering.

We first show the existence of one local transfer $\Pi.$  
If $\pi$ is unramified, then the choice of $\Pi$ is clear; we simply take the irreducible, unramified 
representation determined by the semi-simple conjugacy class in $\GL(N,\C)$ 
generated by the image of the class of $\pi$ under the natural embedding $\iota$ as 
in (\ref{iota}).   We then know, as can be seen directly and is verified in the proof of the global generic transfer 
in the cases we are considering, that we have the requirement of Definition \ref{loc-transfer}, i.e., 
that the local $L$- and $\epsilon$-factors of $\pi$ and $\Pi$ twisted by irreducible, unitary, supercuspidal 
representations $\rho$ of $\GL(m,F)$ for $m$ up to $N-1$ are equal.

Next, assume that $\pi$ is a general (not necessarily unramified) generic supercuspidal representation.   
Since $\pi$ is generic and supercuspidal, by \cite[Proposition 5.1]{shahidi1990-annals}, 
there exist a number field 
$k,$ a non-archimedean place $v_0$ of $k,$ and a globally generic cuspidal 
automorphic representation $\sigma$ of $\G(\mathbb A_k)$ such that 
$k_{v_0} = F,$ and $\sigma_{v_0} = \pi,$ and for all non-archimedean 
places $v \not= v_0$ of $k$ the local representations $\sigma_v$ is 
unramified.  Here, $\sigma$ is generic with respect to a global generic character 
$\Psi$ whose local component at $v_0$ is out fixed $\psi.$ 
(In the non-split quasi-split cases, we have 
a quadratic extension $K/k$ of number fields, and a place $w$ of $K$ 
lying above the place $v$ of $k$ such that $K_w = E.$)

We recall the globally generic automorphic representation $\sigma$ of $\G(\mathbb A_k)$ 
is known to have a transfer $\Sigma$ to $\GL(N,\mathbb A_k)$  
for each of the groups we are considering, as proved in \cite{duke, manuscripta, ckpss1, ckpss2,kk-imrp,kk-aarhus}.  
To be more precise, $\Sigma_v$ is the transfer of $\sigma_v$ as above for $v$ outside a finite set of 
places and $\sigma_v$ unramified.  In particular the twisted $L$- and $\epsilon$-factors are equal for 
$\Sigma_v$ and $\sigma_v$ for such $v.$  Now, just take $\Pi$ to be the local component of $\Sigma$ at $v_0.$

To show this $\Pi$ satisfies the requirements of Definition \ref{loc-transfer},  we again note that 
if $\rho$ is an irreducible, unitary, supercuspidal representation of $\GL(m,F),$ 
$1 \le m \le N-1,$ then we may again use \cite[Proposition 5.1]{shahidi1990-annals} to embed 
$\rho$ in a global cuspidal representation $\tau'$ of $\GL(m,\mathbb A_k)$ such that 
$\tau'_{v_0} = \rho$ and $\tau'_v$ is unramified for all non-archimedean 
$v \not= v_0.$

Let $S$ be a finite set of non-archimedean places of $k$ such that $\sigma_v$ is unramified 
for $v \not\in S$ and let $S' = S - \{v_0\}.$  
Choose an idele class character $\eta$ such that $\eta_{v_0}$ is trivial and $\eta_v$ is highly 
ramified at $v \in S'.$  By a general result, usually referred to as stability of $\gamma$-factors, 
and used in the establishing the global generic transfer results for each of the groups above, 
for $v \in S'$ we have  
\begin{equation}\label{highly-ram} 
\gamma(s, \sigma_v \times (\tau'_v \otimes \eta_v), \Psi_v) 
= 
\gamma(s, \Sigma_v \times (\tau'_v \otimes \eta_v), \Psi_v). 
\end{equation}

Let $\tau = \tau' \otimes \eta.$  Since $\eta_{v_0}$ is trivial, we have 
$\tau_{v_0} = \tau'_{v_0} =\rho.$  On the other hand, we have the global functional equations 
\begin{equation}
L(s, \sigma \times \tau) = \epsilon(s, \sigma \times \tau) 
L(1-s, \w{\sigma} \times \w{\tau})
\end{equation} 
and 
\begin{equation}
L(s, \Sigma \times \tau) = \epsilon(s, \Sigma \times \tau) 
L(1-s, \w{\Sigma} \times \w{\tau}), 
\end{equation} 
which we can rewrite as 
\begin{equation}\label{G-eq}
\gamma(s,\pi \times \rho, \psi) = \left( \prod_{v \in S'} \gamma(s, \sigma_v \times \tau_v, \Psi_v)^{-1} \right) 
\frac{L^S(s,\sigma \times \tau)}{\epsilon^S(s,\sigma \times \tau,\Psi) L^S(1-s, \w{\sigma} \times \w{\tau})} 
\end{equation}
and 
\begin{equation}\label{GL-eq}
\gamma(s,\Pi \times \rho, \psi) = \left( \prod_{v \in S'} \gamma(s, \Sigma_v \times \tau_v, \Psi_v)^{-1} \right) 
\frac{L^S(s,\Sigma \times \tau)}{\epsilon^S(s,\Sigma \times \tau,\Psi) L^S(1-s, \w{\Sigma} \times \w{\tau})}. 
\end{equation}
Now, the fractions on the right hand sides of (\ref{G-eq}) and (\ref{GL-eq}) are equal by the unramified case 
mentioned above and the two products on the right hand sides are also equal as in (highly-ram) above. 
Hence, the left hand sides of (\ref{G-eq}) and (\ref{GL-eq}) must be equal, which means that $\Pi$ is 
indeed a local transfer. 

The uniqueness of $\Pi$ follows from the ``local converse theorem for $\GL(N)$'' which means that 
an irreducible, generic representation of $\GL(N,F)$ (or $\GL(N,E)$ as the case may be) is 
uniquely determined by its $\gamma$-factors twisted by supercuspidal representations of all 
smaller rank general linear groups (cf. Remark after the Corollary of \cite[Theorem 1.1]{henniart}).

It remains to show that $\Pi$ is of the form stated in the theorem.  Being a local component of a globally 
generic automorphic representation, $\Pi$ is a generic, unitary irreducible representation of 
$\GL(N,F)$ or $\GL(N,E)$ as the case may be.  By classification of unitary generic representations of 
the general linear groups \cite{tadic-1986} we have 
\begin{equation}\label{delta-pi}
\Pi = \operatorname{Ind} \left( 
\delta_1 \nu^{r_1} \otimes \cdots \delta_k \nu^{r_k}
\otimes 
\Pi_1 \otimes \cdots \otimes \Pi_d 
\otimes 
\delta_1 \nu^{-r_k}\otimes \cdots \delta_k \nu^{-r_1} 
 \right), 
\end{equation}
where each $\delta_j$ and each $\Pi_i$ is a discrete series representation and $0 < r_k \le \cdots \le r_1 < \frac{1}{2}.$  
Here, $\nu(\cdot) = |\det(\cdot)|.$

Recall that 
\begin{equation}\label{rho-twist}
\gamma(s,\pi\times\rho,\psi) = \gamma(s,\Pi \times \rho, \psi) 
\end{equation}  
for any unitary, supercuspidal representation $\rho$ of $\GL(m)$ for $m$ up to $N-1.$   
In fact, multiplicativity of $\gamma$-factors implies that (\ref{rho-twist}) holds for 
$\rho$ discrete series as well.  To see this,  note that if $\rho$ is discrete series, then it can 
be realized as the irreducible quotient of an induced representation 
\begin{equation} 
\operatorname{Ind}\left(\rho_0 \nu^{-\frac{t-1}{2}}\otimes\cdots\otimes \rho_0 \nu^{\frac{t-1}{2}}\right),  
\end{equation}
where $\rho_0$ is unitary supercuspidal and $t$ is a positive integer. Then, 
\begin{eqnarray*}
\gamma(s,\pi \times \rho, \psi) 
&=& 
\prod_{j=1}^{t-1} \gamma(s+\frac{t-1}{2}-j, \pi\times\rho_0, \psi) \\ 
&=& 
\prod_{j=1}^{t-1} \gamma(s+\frac{t-1}{2}-j, \Pi\times\rho_0, \psi) \\ 
&=& \gamma(s,\Pi\times\rho,\psi),  
\end{eqnarray*}
i.e., (\ref{rho-twist}) holds with $\rho$ any discrete series representation of $\GL(m)$ for $m$ up to $N-1.$

Now, up to a monomial factor coming from the $\epsilon$-factors, the numerator of $\gamma(s, \Pi \times \rho, \psi)$ 
is given by 
\begin{equation}
\left( 
\prod_{j=1}^k L(s+r_j, \delta_j \times \rho) L(s-r_j, \delta_j \times \rho) 
\prod_{i=1}^d L(s, \Pi_i \times \rho)
\right)^{-1}.
\end{equation}
Since neither of $L(s,\delta_j \times \rho)$ or $L(s,\Pi_i \times \rho)$ has a pole in $\Re(s) > 0$ this numerator can have 
zeros only in $\Re(s) < \frac{1}{2}.$ 

Similarly, the denominator of $\gamma(s, \Pi \times \rho, \psi)$ is the polynomial 
\begin{equation}
\left( 
\prod_{j=1}^k L(1-s-r_j, \w{\delta}_j \times \w{\rho}) L(1-s+r_j, \w{\delta}_j \times \w{\rho}) 
\prod_{i=1}^d L(1-s, \w{\Pi}_i \times \w{\rho})
\right)^{-1}, 
\end{equation}
which can only have zeros in the region $\Re(s) > \frac{1}{2}.$

Hence, the numerator and denominator in the factorization coming from the multiplicativity of the $\gamma$-factor 
have not common zeros and, consequently, we conclude from the equality of the $\gamma$-factors that 
\begin{equation}
L(s, \pi \times \rho) = 
\prod_{j=1}^k L(s+r_j, \delta_j \times \rho) L(s-r_j, \delta_j \times \rho) 
\prod_{i=1}^d L(s, \Pi_i \times \rho). 
\end{equation}
On the other hand, we know \cite{jpss-1983} that the same expression gives $L(s,\Pi \times \rho).$  Therefore, 
\begin{equation}\label{L-twist}
L(s,\pi \times \rho) = L(s, \Pi \times \rho), 
\end{equation}
with $\rho$ discrete series.

Fix $1 \le i \le k.$  We apply (\ref{L-twist}) with $\rho = \w{\delta}_j.$  Since $\delta_j$ and $\pi$ are both tempered we 
know that $L(s,\pi \times \w{\delta}_j)$ is holomorphic for $\Re(s) > 0. $  In general, this is the third author's 
Tempered $L$-function Conjecture \cite[Conj. 7.1]{shahidi1990-annals}.  Many cases of this conjecture were proved by several 
authors \cite{As, cas-sha, HO, kim2005-cjm, KKi, muic-shahidi} and a proof in the general case has now appeared in 
\cite{HO}.  
On the other hand, we have 
\begin{equation}
L(s, \Pi \times \w{\delta}_j) = 
\prod_{j=1}^k L(s+r_j, \delta_j \times \w{\delta}_i) L(s-r_j, \delta_j \times \w{\delta}_i) 
\prod_{j=1}^d L(s, \Pi_j \times \w{\delta}_i). 
\end{equation}
The term $L(s-r_i, \delta_i \times \w{\delta}_i)$ produces a pole at $s=r_i$ which results in a pole of 
$L(s,\Pi \times \w{\delta}_i)$ at $s=r_i > 0$ as the local $L$-factors are never zero.  This is a contradiction unless 
$k=0,$ i.e., there are no $\delta_i$'s in (\ref{delta-pi}).   Hence, 
\begin{equation}
\Pi = \operatorname{Ind}\left( \Pi_1 \otimes \cdots \Pi_d \right)  
\end{equation}
is a full induced representation from unitary discrete series representations $\Pi_i.$  In particular, $\Pi$ is tempered.

In fact, we claim that each $\Pi_i$ is unitary supercuspidal.  
To see this, we can again realize the discrete series representation $\Pi_i$ 
as the irreducible quotient of the induced representation 
\begin{equation}
 \operatorname{Ind}\left(\Pi^0_i \nu^{-\frac{t_i-1}{2}} \otimes \cdots \otimes \Pi^0_i \nu^{\frac{t_i-1}{2}} \right)
\end{equation}
associated with the segment $[\Pi^0_i \nu^{-\frac{t_i-1}{2}},\Pi^0_i \nu^{\frac{t_i-1}{2}} ]$ where $\Pi^0_i$ is unitary 
supercuspidal and $t_i$ is a positive integer \cite{bz, zelevinsky-1980}.   Applying (\ref{L-twist}) again with 
$\rho = \w{\Pi}_i$ we have 
\begin{equation}
L(s, \pi \times \w{\Pi}_i) = L(s, \Pi \times \w{\Pi}_i).  
\end{equation}
Let us calculate both sides of this equality.   On the right hand side we have 
\begin{equation}
L(s,\Pi \times \w{\Pi}_i) = \prod_{j=1}^d L(s,\Pi_j \times \w{\Pi}_i) 
\end{equation}
and 
\begin{equation}
L(s,\Pi_i \times \w{\Pi}_i) = \prod_{k=0}^{t_i-1} L(s+k,\Pi^0_i \times \w{\Pi}^0_i).   
\end{equation} 
(The last equation is verified, for example, in \cite[p. 575]{kim2005-cjm}.) 
The local $L$-function $L(s, \Pi^0_i \times \w{\Pi}^0_i)$ has a pole at $s=0$ so that 
$L(s+t_i-1, \Pi^0_i \times \w{\Pi}^0_i)$ has a pole at $s = -(t_i-1).$  Since local $L$-functions are 
never zero, this pole persists to give a pole of $L(s,\Pi \times \w{\Pi}_i)$ at $s=1-t_i.$

As for the left hand side, from \cite{shahidi1990-annals} we have 
\begin{equation}
L(s,\pi\times\w{\Pi}_i) = L(s+\frac{t_i-1}{2}, \pi \times \w{\Pi}_i) 
\end{equation}
since $\pi$ is supercuspidal.   Since $L(s,\pi\times\w{\Pi}_i)$ can have poles only for 
$\Re(s) = 0,$ we see that $L(s,\pi\times\w{\Pi}_i)$ can only have poles on the line 
$\Re(s) = -(t_i-1)/2.$   These locations of poles are inconsistent unless $t_i=1,$ 
i.e., $\Pi_i=\Pi^0_i$ is supercuspidal, as we desire.

Finally, we show that the conditions in terms of the $L$-functions in the statement of the theorem hold.  
Consider the case of the general spin groups first.   Let $\omega = \omega_\pi$ denote 
the central character of $\pi.$   Consider the equality  
\begin{equation}\label{piPi} 
L(s, \pi \times \w{\Pi}_i) = L(s, \Pi \times \w{\Pi}_i). 
\end{equation} 
The right hand side has a pole at $s=0$ as before.  For the left hand side to have a pole at 
$s=0,$ from the general properties of these local $L$-functions (cf. \cite[Cor. 7.6]{shahidi1990-annals}, for example) 
we must have $\Pi \otimes \pi \cong w_0(\Pi \otimes \pi)$ as representations of $\GL \times \G.$  
By Lemma \ref{weyl} below we this implies that 
\begin{equation}
\Pi_i \cong \w{\Pi}_i \otimes \omega.
\end{equation} 
Moreover, the order of pole at $s=0$ on the left hand side of (\ref{piPi}) is one while the order of 
the pole on the right hand side is the number of $j$ such that $\Pi_j \cong \Pi_i.$  Hence, 
$\Pi_i \not\cong \Pi_j$ if $i \not= j.$

Furthermore, assuming that we are in the odd general spin group case, \cite{shahidi1990-annals} 
implies that the product 
\begin{equation}
L(s, \pi \times \Pi_i) L(s, \Pi_i, \sym^2 \otimes \omega^{-1}) 
\end{equation} 
has a simple pole at $s=0$ in this situation.  This pole is already accounted for by the pole 
at $s=0$ of $L(s, \pi \times \Pi_i).$   Therefore, $L(s, \Pi_i, \sym^2 \otimes \omega^{-1})$ has no 
pole at $s=0.$  However, 
\begin{equation}
L(s, \Pi_i \times \w{\Pi}_i) = L(s, \Pi_i \times \Pi_i \otimes \omega^{-1}) = 
L(s, \Pi_i, \sym^2 \otimes \omega^{-1}) L(s, \Pi_i, \wedge^2 \otimes \omega^{-1}),  
\end{equation}
which implies that $L(s, \Pi_i, \wedge^2 \otimes \omega^{-1})$ has a pole at $s=0$ (which can only happen 
if $N_i$ is even).  If we are in the even general spin groups, the same argument works with the roles 
of $\sym^2$ and $\wedge^2$ switched.

We end the proof by mentioning that a similar argument, with minor modifications replacing $\omega$-self-dual 
with self-dual or conjugate-self-dual as appropriate, establishes the $L$-function condition for 
the remaining groups in the statement of the theorem.   We will not repeat those arguments as they 
are similar and have already appeared in the literature. For orthogonal and symplectic 
groups, this is done in \cite[p. 203]{ckpss2}.  For unitary groups, it is verified in \cite{kk-imrp} and 
\cite{kk-aarhus}.  
\end{proof}

We should note here that the conditions that the transfers $\Pi$ need to satisfy in the theorem 
above in fact determine the image of the transfer.  In other words, every $\Pi$ satisfying these 
conditions is the transfer of an irreducible, generic, supercuspidal $\pi$ from the  
appropriate group $\G$ to $\GL(N).$  For this one needs the ``local descent'' for all the groups 
we are considering.  For symplectic and orthogonal groups, as well as the unitary groups, these 
facts have already been established \cite{soudry-asterisque, grs} and for the general spin groups this is a 
work in progress of Jing Feng Lau.

\section{Reducibility of Local Representations} \label{loc-red}

As an application of our results on the generic local transfer in Section \ref{sc-transfer}, we now 
give some uniform results on reducibility of local induced representations. 

\subsection{Reducibility on the Unitary Axis} \label{sec-unit}

To state our main uniform results on irreducibility, we first define the orthogonal/symplectic representations of 
general linear groups. 

Let $F$ and $\psi$ be as before. 
Let $\eta$ be a character of $F^\times$ and let $\sigma$ be an irreducible 
supercuspidal representation of $\GL(m,F).$  Let 
\begin{equation} 
\phi : W_F \longrightarrow \GL(m,\C) 
\end{equation}
be the parameter of $\sigma$ and set 
\begin{eqnarray}
L(s, \sigma, \sym^2\otimes \eta) &:=& L(s, \sym^2 \phi \cdot \eta), \\
L(s, \sigma, \wedge^2\otimes \eta) &:=& L(s, \wedge^2 \phi \cdot \eta).  
\end{eqnarray} 
When $\eta = 1,$ these reduce to the usual untwisted $L$-factors. 

Similarly, with $E/F$ as before, let $\sigma$ be an irreducible 
supercuspidal representation of $\GL(m,E).$ Let 
\begin{equation} 
\phi : W_E \longrightarrow \GL(m,\C) 
\end{equation}
be the parameter of $\sigma$ and set 
\begin{eqnarray} 
L(s, \sigma, r_A) &:=& L(s, r_A\circ \phi), \\
L(s, \sigma, r_A \otimes \delta_{E/F}) &:=& L(s, r_A\circ \phi \cdot \delta_{E/F}). 
\end{eqnarray} 
Here, $r_A$ is denote the Asai representation.  For details about the Asai representation and 
the corresponding $L$-function we refer to \cite[\S 5 and \S 8]{kk-imrp}. 

The representation $\sigma$ of $\GL(m,F)$ is said to be $\eta$-self-dual 
if it satisfies 
\begin{equation} 
\w{\sigma} \cong \sigma \otimes \eta. 
\end{equation}  
If $\sigma$ is $\eta$-self-dual, then we have 
\begin{eqnarray} \label{omega-self-dual}
L(s, \sigma \times \w{\sigma}) &=& L(s, \sigma \times \sigma \otimes \eta) \nonumber \\ 
&=& L(s, \sigma, \wedge^2\otimes \eta)  \cdot L(s, \sigma, \sym^2\otimes \eta) 
\end{eqnarray}
and exactly one of the two local $L$-function on the right hand side of (\ref{omega-self-dual}) 
has a pole at $s=0.$  Conversely, if one of the $L$-functions on the right hand side of 
(\ref{omega-self-dual}) has a pole at $s=0,$ then $\sigma$ is $\eta$-self-dual. 
Again, when $\eta=1$ the representation $\sigma$ is said to be self-dual and the (untwisted) 
exterior/symmetric square $L$-functions replace the twisted ones above. 

Analogously, a representation $\sigma$ of $\GL(m,E)$ is said to be conjugate-self-dual if it satisfies 
\begin{equation} 
\w{\sigma} \cong \overline{\sigma}. 
\end{equation}  
We recall that $x \mapsto \overline{x}$ denotes the action of $\operatorname{Gal}(E/F)$ on the 
representation $\sigma$ of $\GL(m,E)$ and $\overline{\sigma}$ denotes the corresponding action 
on $\sigma.$  For $\sigma$ conjugate-self-dual, we have 
\begin{eqnarray} \label{conj-self-dual}
L(s, \sigma \times \w{\sigma}) &=& L(s, \sigma \times \overline{\sigma}) \nonumber \\ 
&=& L(s, \sigma, r_A)  \cdot L(s, \sigma, r_A \otimes \delta_{E/F}),  
\end{eqnarray}
where the local $L$-functions on the right hand side are, as before, the Asai $L$-function and its 
twist by the quadratic character $\delta_{E/F}.$  Again, the $L$-function on the 
left hand side of (\ref{conj-self-dual}) has a pole at $s=0$ which implies that exactly one of 
those on the right hand side of (\ref{conj-self-dual}) has a pole at $s=0.$ 

\begin{defi} \label{os}
An irreducible, unitary, supercuspidal representation 
$\sigma$ of $\GL(m,F)$ is said to be $\eta$-symplectic if 
$L(s, \sigma, \wedge^2\otimes \eta)$ has a pole at $s=0,$ and it is 
said to be $\eta$-orthogonal if $L(s, \sigma, \sym^2\otimes \eta)$ has a pole at $s=0.$ 
When $\eta=1,$ then $\sigma$ is simply called symplectic or orthogonal. 

Similarly, an irreducible, unitary, supercuspidal representation 
$\sigma$ of $\GL(m,E)$ is said to be conjugate-orthogonal if 
$L(s, \sigma, r_A)$ has a pole at $s=0,$ and it is said to be conjugate-symplectic if 
$L(s, \sigma, r_A \otimes \delta_{E/F})$ has a pole at $s=0.$ 

\end{defi} 

Note that if a representation $\sigma$ of $\GL(m,F)$ is either $\eta$-symplectic or 
$\eta$-orthogonal, then (\ref{omega-self-dual}) 
implies that $\sigma$ is $\eta$-self-dual. Moreover, in the $\eta$-symplectic case, $m$ must 
be even.  In the following theorem we use these notions with $\eta = \omega^{-1},$ where 
$\omega = \omega_\pi$ denotes the central character of a representation of $\G(F).$ 
Similarly, if a representation $\sigma$ of $\GL(m,E)$ is either conjugate-symplectic 
or conjugate-orthogonal, then (\ref{conj-self-dual}) 
implies that $\sigma$ is conjugate-self-dual.

Let $\G=\G(n)$ be as in Section \ref{notation}, i.e., 
$\G$ is one of the groups 
$\SO(2n+1)$, $\Sp(2n)$, $\SO(2n)$, $\GSpin(2n+1)$, $\GSpin(2n),$ 
or 
$\SO_{E/F}(2n), \GSpin_{E/F}(2n),$  
or 
$\Uef(2n+1), \Uef(2n),$ 
where $E/F$ is a 
quadratic extension. 
Let $\widehat{\G}$ denote the connected component of its Langlands $L$-group ${}^L\G.$  
We define the type of $\widehat{\G}$ as follows: 

\[ 
\renewcommand{\arraystretch}{1.5}
\begin{array}{|c|c|c|}
\hline
\G & \widehat{\G} & \mbox{type of } \widehat{\G} \\
\hline
\SO(2n+1) &  \Sp(2n,\C) & \mbox{symplectic} \\
\Sp(2n) &  \SO(2n+1,\C) & \mbox{orthogonal} \\
\SO(2n) &  \SO(2n,\C) & \mbox{orthogonal} \\
\GSpin(2n+1) &  \GSp(2n,\C) & \mbox{symplectic} \\
\GSpin(2n) &  \GSO(2n,\C) & \mbox{orthogonal} \\
\SO_{E/F}(2n)& \SO(2n,\C) & \mbox{orthogonal} \\ 
\GSpin_{E/F}(2n) & \GSO(2n,\C) & \mbox{orthogonal} \\
\Uef(2n) 
& 
\GL(2n,\C) 
& 
\mbox{(conjugate) symplectic}  
\\
\Uef(2n+1) 
& 
\GL(2n+1,\C) 
& 
\mbox{(conjugate) orthogonal}  
\\
\hline
\end{array}
\]

\begin{thm} \label{reduce-unitary}
Let $n \ge 0$ and $m \ge 1$ be integers. 
Let $\G=\G(n)$ and $E/F$ be as before. 
Let $\sigma$ be an irreducible, unitary, supercuspidal representation of $\GL(m,F)$ 
if $\G$ is orthogonal or symplectic type, or of $\GL(m,E)$ when $\G$ is unitary.    
Let $\pi$ be an irreducible, generic, unitary, supercuspidal representation of $\G(n,F)$ with 
central character $\omega = \omega_\pi.$   
Consider $ \tau = \sigma \otimes \pi$, an irreducible, generic, unitary, supercuspidal 
representation of $\M(F),$ where $\M=\GL(m) \times \G(n)$ or 
$\M=\operatorname{Res}_{E/F} \GL(m) \times \G(n),$ as appropriate, 
is a standard Levi subgroup of a 
maximal parabolic  subgroup $\P$ in the larger quasi-split group $\G(m+n).$ 
Let $I(\tau) = I(\sigma \otimes \pi)$ be the induced representation of $\G(m+n,F)$ as in (\ref{ind-0}).

If $\P$ is not self-associate (cf. Remark \ref{non-self-assoc}), then $I(\tau)$ is irreducible .  

If $\P$ is self-associate, then for $\G$ orthogonal, symplectic or a general spin group, we have the following statements: 

\begin{itemize}

\item[(a)] If $\sigma$ is not $\omega^{-1}$-self-dual, then $I(\sigma\otimes\pi)$ is irreducible. 

\item[(b)] If $\sigma$ is $\omega^{-1}$-self-dual and of the type opposite to $\widehat \G$, 
then $I(\sigma\otimes\pi)$ is irreducible. 

\item[(c)] If $\sigma$ is $\omega^{-1}$-self-dual and of the same type as $\widehat \G$, 
then $I(\sigma\otimes\pi)$ is irreducible if and only if 
$\sigma$ appears as a component of the transfer of $\pi$ to $\GL(N,F)$ as in Theorem \ref{local-transfer-sc}. 
In particular, if $m > N,$ then $I(\sigma \otimes \pi)$ is always reducible. 

\end{itemize}

Analogously, for $\G$ unitary, we have the following statements: 

\begin{itemize}

\item[(d)] If $\sigma$ is not conjugate-self-dual, then $I(\sigma\otimes\pi)$ is irreducible. 

\item[(e)] If $\sigma$ is conjugate-self-dual and of the type opposite to $\widehat \G$, 
then $I(\sigma\otimes\pi)$ is irreducible. 

\item[(f)] If $\sigma$ is conjugate-self-dual and of the same type as $\widehat \G$, 
then $I(\sigma\otimes\pi)$ is irreducible if and only if 
$\sigma$ appears as a component of the transfer of $\pi$ to $\GL(N,E)$ as in Theorem \ref{local-transfer-sc}. 
In particular, if $m > N,$ then $I(\sigma \otimes \pi)$ is always reducible. 

\end{itemize}

\end{thm}

\begin{proof}
The theorem essentially follows from Corollary \ref{red-cor} combined with the local 
transfer results Theorem \ref{local-transfer-sc} as we now explain. 

If $\P$ is not self-associate, then clearly $w_0(\tau) \not\cong \tau$ and 
$I(\tau)$ is irreducible by Corollary \ref{red-cor}.  

Next, assume that $\P$ is self-associate. We then have 
\begin{equation}
w_0 ( \sigma \otimes \pi) \cong \w{\sigma} \cdot \omega_\pi \otimes \pi.
\end{equation} 
when $\G$ is orthogonal, symplectic, or a general spin group.  Similarly, 
\begin{equation}
w_0 ( \sigma \otimes \pi) \cong \w{\overline{\sigma}} \otimes \pi.
\end{equation} 
when $\G$ is unitary (cf. Lemma \ref{weyl}). 

Therefore, if $\G$ is orthogonal or symplectic type and $\sigma$ is not $\omega^{-1}$-self-dual, 
or if $\G$ is unitary and $\sigma$ is not conjugate-self-dual, then 
$w_0 ( \sigma \otimes \pi) \not\cong \sigma \otimes \pi,$ and part (a) of 
Corollary \ref{red-cor} 
implies that $I(\sigma\otimes\pi)$ is irreducible.  This proves statements (a) and (d).

Now, let $\G=\G(n)=\SO(2n+1)$ or $\GSpin(2n+1)$ and let $\M=\GL(m) \times \G(n)$ 
as a Levi in $\G(m+n).$  Then $\widehat\G$ is of symplectic type. 
Consider the case of $n \ge 1$ (non-Siegel Levi) first. 
Assume that $\sigma$ is $\omega^{-1}$-self-dual.  Then $\sigma$ is either 
$\omega^{-1}$-orthogonal or $\omega^{-1}$-symplectic (cf. Definition \ref{os}).  
If $\sigma$ is of 
the type opposite to $\widehat \G,$ then $\sigma$ is $\omega^{-1}$-orthogonal, which 
means that the local $L$-function $L(s,\sigma, \sym^2\otimes\omega^{-1})$ has a pole 
at $s=0.$  This local $L$-function is the second $L$-function appearing in 
Theorem \ref{const-term} or part (b) of Corollary \ref{red-cor}.  On the other hand, 
if $\sigma$ is of the same type as $\widehat \G,$ then it is $\omega^{-1}$-symplectic 
and, hence, $L(s,\sigma,\wedge^2\otimes\omega^{-1})$ has a pole at $s=0$ and 
$L(s,\sigma, \sym^2\otimes\omega^{-1})$ does not.   Now, the  other (first) $L$-function in part (b) of Corollary \ref{red-cor} 
would have a pole at $s=0$ if and only if $\sigma$ appears as a component in the transfer $\Pi$ of $\pi.$ 
To see this, note that the other $L$-function is 
\begin{equation} 
L(s,\w{\sigma} \times \pi) = L(s,\w{\sigma} \times \Pi) 
= \prod_{i=1}^d L(s, \w{\sigma} \times \Pi_i), 
\end{equation}
where $\Pi_i$'s are the components of the transfer $\Pi$ of $\pi$ as in Theorem \ref{local-transfer-sc}.

If $n=0,$ the group $\G(0)$ is either trivial in which case $\pi$ is trivial, or isomorphic to $\GL(1)$ in which 
case $\pi = \omega_\pi$ is just a character.  This is the Siegel Levi case and in this case only one 
$L$-function, $L(s,\sigma, \sym^2\otimes\omega^{-1}),$ appears in Corollary \ref{red-cor}.  
The above argument still holds in the following sense.  If $\sigma$ is $\omega^{-1}$-orthogonal, 
then the first (and only) local $L$-function in Corollary \ref{red-cor} has a pole at $s=0$ and 
$I(\sigma\otimes\pi)$ is irreducible.  If $\sigma$ is $\omega^{-1}$-symplectic, then $I(\sigma\otimes\pi)$ 
is reducible.  Neither does any $L$-function in Corollary \ref{red-cor} have a pole at $s=0,$ nor does 
$\sigma$ appear as a component of transfer $\Pi$ of $\pi.$ 
This proves (b) and (c) for $\G=\G(n)=\SO(2n+1)$ or $\GSpin(2n+1).$

Next, let $\G=\G(n) = \SO(2n),$ $\GSpin(2n)$ or their quasi-split forms.  A similar 
argument as above again holds, except that $\widehat \G$ is now of orthogonal type and if 
$\sigma$ is of the type opposite to $\widehat \G,$ then it is $\omega^{-1}$-symplectic, which 
means that the local $L$-function $L(s,\sigma,\wedge^2\otimes\omega^{-1})$ has a pole 
at $s=0.$  Now, this is the second $L$-function appearing in part (b) of Corollary \ref{red-cor}.  
And if $\sigma$ is of the same type at $\widehat \G,$ then it is $\sigma^{-1}$-orthogonal 
and, hence, $L(s,\sigma,\sym^2\otimes\omega^{-1})$ has a pole at $s=0$ and 
$L(s,\sigma,\wedge^2\otimes\omega^{-1})$ does not.  Now, in a similar way, the other $L$-function in 
part (b) of Corollary \ref{red-cor} would have a pole at $s=0$ if and only if 
$\sigma$ appears as a component of the transfer $\Pi$ of $\pi.$  When $n=0$ a similar situation 
occurs with one one local $L$-function appearing again. 

When $\G=\G(n) = \Sp(2n)$, the above paragraph holds again.  The difference is just that the 
transfer $\Pi$ is a representation of $\GL(2n+1,F).$  When $n \ge 1$ there are two $L$-functions, 
namely, $L(s, \sigma \times \pi)$ and $L(s, \sigma, \wedge^2).$  When $n=0,$ there are 
actually again two $L$-functions appearing, namely, $L(s,\sigma)$ (the standard $L$-function) which 
does not produce any poles at $s=0,$ and $L(s,\sigma, \wedge^2)$ which behaves the same way as 
above.  

Hence, we have proved parts (b) and (c) for $\G= \Sp(2n), \SO(2n), \GSpin(2n)$ and their 
quasi-split forms. 

Finally, let $\G=\G(n) = \Uef(2n)$ or $\Uef(2n+1)$ and assume that the representation 
$\sigma$ of $\GL(m,E)$ is conjugate-self-dual.  Then $\sigma$ is either 
conjugate-orthogonal or conjugate-symplectic (cf. Definition \ref{os}).  

Consider $\G=\Uef(2n)$ first. Now, $\widehat \G$ is (conjugate) symplectic. 
If $\sigma$ is of type opposite to $\widehat\G,$ then it is conjugate-orthogonal, 
which means that the local Asai $L$-function 
$L(s,\sigma, r_A)$ has a pole at $s=0.$  
This $L$-function is the second $L$-function appearing in 
Theorem \ref{const-term} or part (b) of Corollary \ref{red-cor}.  
On the other hand, 
if $\sigma$ is of the same type as $\widehat \G,$ then it is conjugate-symplectic 
and, hence, $L(s,\sigma, r_A \otimes\delta_{E/F})$ has a pole at $s=0$ and 
$L(s,\sigma, r_A)$ does not.   Now, in a similar way as above, the other 
$L$-function in part (b) of Corollary \ref{red-cor} would have a pole at $s=0$ 
if and only if $\sigma$ appears as a component of the transfer $\Pi$ of $\pi.$ 

The argument for $\G=\Uef(2n+1)$ is exactly the same with the words 
(conjugate) symplectic and (conjugate) orthogonal switched.  

Therefore, we have also proved (e) and (f) for $\G=\Uef(2n)$ and $\Uef(2n+1),$ 
which finishes the proof of the theorem. 
\end{proof}

\begin{rem}
It is worth pointing out that in the case of $\G=\GSp(2n),$ if the representation $\sigma$ is 
$\omega^{-1}$-self-dual, then $\omega=1$ (cf. \cite[p. 286]{shahidi91pm}).  The case of 
non-trivial $\omega$ may occur for general spin groups. 
\end{rem}

\begin{lem} \label{weyl}
Let $m$ and $n$ be non-negative integers and let $\G = \G(m+n)$ and $E/F$ be as before. 
Let $\theta = \Delta - \{\alpha \},$ where $\Delta$ denotes the set of simple roots of 
$\G$ and $\alpha$ is a fixed simple root.  
Consider the standard maximal parabolic subgroup $\P= \P_\theta = \M \N$ with 
the Levi $\M \cong \GL(m) \times \G(n)$ if $\G$ is one of the non-unitary groups 
we are considering, or 
$\M \cong \operatorname{Res}_{E/F} \GL(m) \times \G(n)$ if $\G$ is unitary. 
Let $w_0$ be the unique element in the Weyl group of $\G$ such that $w_0(\theta) \subset \Delta$ 
and $w_0 (\alpha) < 0.$  
We assume that $\P$ is self-associate, i.e., $w_0(\theta) = \theta$ (cf.  Remark \ref{non-self-assoc}).  

Let $\sigma$ be a representation of $\GL(m,F),$ or of $\GL(m,E)$ when $\G$ is unitary, and 
let $\pi$ 
be a representation of $\G(n,F).$ 

Then, 
\[
w_0 (\sigma \otimes \pi) \cong \w{\sigma} \otimes \pi. 
\] 
when $\G = \SO(2n+1), \SO(2n),$ or $\SO_{E/F}(2n),$ 
and  
\[
w_0 (\sigma \otimes \pi) \cong \w{\sigma} \cdot \omega_\pi \otimes \pi. 
\] 
when $\G = \GSpin(2n+1), \GSpin(2n),$ or $\GSpin_{E/F}(2n).$ 
 
For $\G=\Uef(2n)$ or $\Uef(2n+1),$ we have 
\[
w_0 (\sigma \otimes \pi) \cong \w{\overline{\sigma}} \otimes \pi. 
\] 
\end{lem}

\begin{proof}
One verifies this lemma by considering the effect of conjugation by the Weyl group element $w_0$ on 
an element of the Levi $\M.$   For special orthogonal, symplectic, or unitary groups, we can do this by 
a standard matrix calculation, noting that the action of $w_0$ is to simply switch the upper left $m \times m$ 
block with the lower right block of the same size in the usual matrix representation of these groups.  

For the general spin groups essentially the same observation works, except that one expresses it in 
terms of root data due to lack of a convenient matrix realization and follows the action of the Weyl group 
element $w_0.$  Let us give some details for this case. 

Consider $\G = \G(m+n) = \GSpin(2m+2n+1).$ 
Using the Bourbaki notation, a detailed description of the root data for $\G$ is given in \cite[\S 1]{manuscripta} 
which we use below.  
Let 
\begin{equation} 
X = \Z \langle e_0, e_1, \dots, e_{m+n} \rangle
\end{equation} 
and 
\begin{equation} 
X^\vee = \Z \langle e^*_0, e^*_1, \dots, e^*_{m+n} \rangle
\end{equation} 
denote the character and cocharacter lattices of $\G,$ respectively, with the standard $\Z$-pairing.  With the simple 
roots $\Delta = \{ \alpha_1, \dots, \alpha_{m+n} \}$ and the simple coroots 
$\Delta^\vee = \{ \alpha^\vee_1, \dots, \alpha^\vee_{m+n} \}$ defined as in \cite[\S 1]{manuscripta}, we have 
$\theta = \Delta - \{ \alpha_{m} \}$ with $w_0$ as in Section \ref{notation}.   
Let $\M = \M_\theta$ be the maximal standard Levi subgroup corresponding to $\theta.$Then 
\begin{equation}\label{levi-iso}
\M \cong \GL(m) \times \G(n)
\end{equation} 
with 
$X = X_1 \oplus X_2$ and $X^\vee = X^\vee_1 \oplus X^\vee_2,$ where 
\begin{equation}
X_1 = \Z \langle e_1, \dots, e_{m} \rangle \quad 
X_2 = \Z \langle e_0, e_{m+1}, \dots, e_{m+n} \rangle 
\end{equation} 
and 
\begin{equation}
X^\vee_1 = \Z \langle e^\vee_1, \dots, e^\vee_{m} \rangle \quad 
X^\vee_2 = \Z \langle e^\vee_0, e^\vee_{m+1}, \dots, e^\vee_{m+n} \rangle. 
\end{equation} 
Now if we translate the action of $w_0$ on the root data from $\M$ to $\GL(m) \times \G(n)$ via 
the isomorphism (\ref{levi-iso}), we can conclude that for $m=(A,g)$ with $A \in \GL(m)$ and 
$g \in \G(n)$ we have 
\begin{equation}\label{w0-mu}
w_0(m) = ( \mu \cdot {}^tA^{-1} , g), 
\end{equation}
where $\mu = e_0(g)$ is the ``similitude character''.  
This proves the statement of the lemma in this case. 

The case of even general spin groups is similar.  However,  in the even case (\ref{w0-mu}) 
holds provided that we are in the self-associate case (cf. Remark \ref{non-self-assoc}).  This proves 
the lemma.  
\end{proof}

\subsection{Reducibility off the Unitary Axis} \label{sec-off-unit}
Theorem \ref{reduce-unitary} determines the reducibility of representations of classical groups induced 
from unitary, generic, supercuspidal representation of a maximal Levi in a satisfactory way.  
The analogous question for when the inducing representation is non-unitary is fortunately 
reduced to the unitary case thanks to the following rather general theorem, a  well-known 
result in the Langlans-Shahidi method (cf. \cite[Theorem 8.1]{shahidi1990-annals} or 
\cite[Theorem 5.1]{shahidi91pm}, for example.)

\begin{thm} \label{reduce-off-unitary}
Let $\G,$ $\P=\M\N,$ $\tau=\sigma\otimes\pi$ and $w_0$ be as before.  
Assume that $w_0(\tau) \cong \tau$ and that $I(\tau)$ is 
irreducible.  Let $i=1$ or $2$ be the unique index such that $L(s,\tau,\w{r}_i)$ has a pole at $s=0$ as in 
Corollary \ref{red-cor}.  Then, the induced representation $I(s,\tau)$ of (\ref{ind-s}) is 
\begin{itemize}
\item[(a)] irreducible for $0 < s < 1/i.$ 
\item[(b)] reducible for $s = 1/i.$ 
\item[(c)] irreducible for $s > 1/i.$ 
\end{itemize}
If $w_0(\tau) \cong \tau$ and $I(\tau)$ reduces, 
then $I(s,\tau)$ is irreducible for $s>0.$
\end{thm}

Recall that for the groups we are considering, we always have $i=1$ or $i=2.$  Hence, the 
reducibility point of the induced representation $I(s,\tau)$ is always at either 
$s=1/2$ or $s=1,$ if any, in the region $s > 0.$  In Section \ref{sec-gps} we specify these 
reducibility points for each group individually.

Moreover, we also recall, as one checks easily using the roots of $\G$ in each case, 
that the following equalities are immediate from (\ref{ind-s}): 

\begin{equation} \label{det-s}
I(s, \tau) = \begin{cases}
	I\left(\nu^{s} \sigma \otimes \pi \right) & \mbox{ if $\G(n)$ is of type $A$ (unitary),}  \\ 
	I\left(\nu^{s} \sigma \otimes \pi \right) & \mbox{ if $\G(n)$ is of type $B$ or $D$ and $n \ge 1,$} \\
	I\left(\nu^{s/2} \sigma \otimes \pi \right) & \mbox{ if $\G(n)$ is of type $B$ or $D$ and $n = 0,$} \\
	I\left(\nu^{s} \sigma \otimes \pi\right) & \mbox{ if $\G(n)$ is of type $C.$} 
		\end{cases}
\end{equation}
Here, $\nu=|det|$ denotes the $p$-adic absolute value of the determinant character on $\GL(m,F)$  
(or $\GL(m,E)$ as the case may be).   When we summarize our 
reducibility results for each individual group in Section \ref{sec-gps}, we will state them in terms 
of $\det$ rather than $\w{\alpha}$ in (\ref{ind-s}).

\subsection{Reducibility for Groups of Classical Type} \label{sec-gps}
We now summarize our results on reducibility points of the induced representations 
for each of the groups we consider in this article.  Below $F$ continues to denote 
a non-archimedean local field of characteristic zero and, when appropriate, 
$E/F$ denotes a quadratic extension, as before.

\begin{prop}\label{so-odd} Reducibility for $\SO(2n+1)$. \newline 
Let $m \ge 1$ and let $\sigma$ be an irreducible, unitary, supercuspidal 
representation of $\GL(m, F).$  Let $n \ge 0$ and let $\pi$ be an irreducible, 
generic, unitary, supercuspidal representation of $\SO(2n+1,F).$  Let 
$I(s) = \operatorname{Ind}\left(|\det|^s \sigma \otimes \pi \right)$ denote the parabolically 
induced representation of $\SO(2m+2n+1,F).$  The following hold: 

\begin{itemize}

\item If $\sigma$ is not self-dual, then $I(s)$ is irreducible for $s \ge 0.$ 

\item If $\sigma$ is self-dual and $L(s, \sigma, \sym^2)$ has a pole at $s=0,$ then $I(s)$ 
is irreducible for $ 0 \le s < 1/2,$ reducible for $s=1/2$, and irreducible for $s > 1/2.$ 
(The reducibility point is the same whether $n = 0$ or $n \ge 1.$) 

\item If $\sigma$ is self-dual, $L(s, \sigma, \wedge^2)$ has a pole at $s=0$ and $n=0,$ 
then $I(s)$ is reducible for $s=0$ and irreducible for $s > 0.$ 

\item If $\sigma$ is self-dual, $L(s, \sigma, \wedge^2)$ has a pole at $s=0,$ $n \ge 1,$ and 
$\sigma$ appears as a component of the transfer of $\pi$ to $\GL(2n,F),$  then $I(s)$ 
is irreducible for $0 \le s < 1,$ reducible for $s=1,$ and irreducible for $s >1.$ 

\item If $\sigma$ is self-dual, $L(s, \sigma, \wedge^2)$ has a pole at $s=0,$ $n \ge 1$ and 
$\sigma$ does not appear as a component of the transfer of $\pi,$ then $I(s)$ 
is reducible for $s=0$ and irreducible for $s > 0.$ 

\end{itemize}

\end{prop}

\begin{prop}\label{gspin-odd} Reducibility for $\GSpin(2n+1).$  \newline
Let $m \ge 1$ and let $\sigma$ be an irreducible, unitary, supercuspidal 
representation of $\GL(m, F).$  Let $n \ge 0$ and let $\pi$ be an irreducible, 
generic, unitary, supercuspidal representation of $\GSpin(2n+1,F)$ with central 
character $\omega = \omega_\pi.$   Let 
$I(s) = \operatorname{Ind}\left(|\det|^s \sigma \otimes \pi \right)$ denote the parabolically 
induced representation of $\GSpin(2m+2n+1,F).$  The following hold: 

\begin{itemize}

\item If $\sigma$ is not $\omega^{-1}$-self-dual, then $I(s)$ is irreducible for $s \ge 0.$ 

\item If $\sigma$ is $\omega^{-1}$-self-dual and $L(s, \sigma, \sym^2\otimes\omega^{-1})$ has a pole at $s=0,$ then $I(s)$ 
is irreducible for $ 0 \le s < 1/2,$ reducible for $s=1/2$, and irreducible for $s > 1/2.$ 
(The reducibility point is the same whether $n = 0$ or $n \ge 1.$) 

\item If $\sigma$ is $\omega^{-1}$-self-dual, $L(s, \sigma, \wedge^2\otimes\omega^{-1})$ has a pole at $s=0$ and $n=0,$ 
then $I(s)$ is reducible for $s=0$ and irreducible for $s > 0.$ 

\item If $\sigma$ is $\omega^{-1}$-self-dual, $L(s, \sigma, \wedge^2\otimes\omega^{-1})$ has a pole at $s=0,$ $n \ge 1,$ and 
$\sigma$ appears as a component of the transfer of $\pi$ to $\GL(2n,F),$  then $I(s)$ 
is irreducible for $0 \le s < 1,$ reducible for $s=1,$ and irreducible for $s >1.$ 

\item If $\sigma$ is $\omega^{-1}$-self-dual, $L(s, \sigma, \wedge^2\otimes\omega^{-1})$ has a pole at $s=0,$ $n \ge 1$ and 
$\sigma$ does not appear as a component of the transfer of $\pi,$ then $I(s)$ 
is reducible for $s=0$ and irreducible for $s > 0.$ 

\end{itemize}

\end{prop}

\begin{prop}\label{sp} Reducibility for $\Sp(2n).$ \newline
Let $m \ge 1$ and let $\sigma$ be an irreducible, unitary, supercuspidal 
representation of $\GL(m, F).$  Let $n \ge 0$ and let $\pi$ be an irreducible, 
generic, unitary, supercuspidal representation of $\Sp(2n,F).$  Let 
$I(s) = \operatorname{Ind}\left(|\det|^s \sigma \otimes \pi \right)$ denote the parabolically 
induced representation of $\Sp(2m+2n,F).$  The following hold: 

\begin{itemize}

\item If $\sigma$ is not self-dual, then $I(s)$ is irreducible for $s \ge 0.$ 

\item If $\sigma$ is self-dual and $L(s, \sigma, \wedge^2)$ has a pole at $s=0,$ then $I(s)$ 
is irreducible for $ 0 \le s < 1/2,$ reducible for $s=1/2$, and irreducible for $s > 1/2.$ 
(The reducibility point is the same whether $n = 0$ or $n \ge 1.$) 

\item If $\sigma$ is self-dual, $L(s, \sigma, \sym^2)$ has a pole at $s=0$ and $n=0,$ 
then $I(s)$ is reducible for $s=0$ and irreducible for $s > 0.$ 

\item If $\sigma$ is self-dual, $L(s, \sigma, \wedge^2)$ has a pole at $s=0,$ $n \ge 1,$ and 
$\sigma$ appears as a component of the transfer of $\pi$ to $\GL(2n,F),$  then $I(s)$ 
is irreducible for $0 \le s < 1,$ reducible for $s=1,$ and irreducible for $s >1.$ 

\item If $\sigma$ is self-dual, $L(s, \sigma, \wedge^2)$ has a pole at $s=0,$ $n \ge 1$ and 
$\sigma$ does not appear as a component of the transfer of $\pi,$ then $I(s)$ 
is reducible for $s=0$ and irreducible for $s > 0.$ 

\end{itemize}

\end{prop}

\begin{prop}\label{so-even} Reducibility for $\G(n) = \SO(2n)$ or  $\SO_{E/F}(2n).$ \newline 
Let $m \ge 1$ and let $\sigma$ be an irreducible, unitary, supercuspidal 
representation of $\GL(m, F).$  Let $n \ge 0$ and let $\pi$ be an irreducible, 
generic, unitary, supercuspidal representation of $\G(n,F).$  
Let $I(s) = \operatorname{Ind}\left(|\det|^s \sigma \otimes \pi \right)$ denote the parabolically 
induced representation of $\G(m+n,F).$  

If $n=0$ and $m$ is odd, i.e., non-self-associate parabolic (cf. Remark \ref{non-self-assoc}), 
then $I(s)$ is always irreducible.  Otherwise, the following hold: 

\begin{itemize}

\item If $\sigma$ is not self-dual, then $I(s)$ is irreducible for $s \ge 0.$ 

\item If $\sigma$ is self-dual and $L(s, \sigma, \wedge^2)$ has a pole at $s=0,$ then $I(s)$ 
is irreducible for $ 0 \le s < 1/2,$ reducible for $s=1/2$, and irreducible for $s > 1/2.$ 
(The reducibility point is the same whether $n = 0$ or $n \ge 1.$) 

\item If $\sigma$ is self-dual, $L(s, \sigma, \sym^2)$ has a pole at $s=0$ and $n=0,$ 
then $I(s)$ is reducible for $s=0$ and irreducible for $s > 0.$ 

\item If $\sigma$ is self-dual, $L(s, \sigma, \sym^2)$ has a pole at $s=0,$ $n \ge 1,$ and 
$\sigma$ appears as a component of the transfer of $\pi$ to $\GL(2n,F),$  then $I(s)$ 
is irreducible for $0 \le s < 1,$ reducible for $s=1,$ and irreducible for $s >1.$ 

\item If $\sigma$ is self-dual, $L(s, \sigma, \sym^2)$ has a pole at $s=0,$ $n \ge 1$ and 
$\sigma$ does not appear as a component of the transfer of $\pi,$ then $I(s)$ 
is reducible for $s=0$ and irreducible for $s > 0.$ 

\end{itemize}

\end{prop}

\begin{prop}\label{gspin-even} 
Reducibility for $\G(n) = \GSpin(2n)$ and $\GSpin_{E/F}(2n).$ \newline
Let $m \ge 1$ and let $\sigma$ be an irreducible, unitary, supercuspidal 
representation of $\GL(m, F).$  Let $n \ge 0$ and let $\pi$ be an irreducible, 
generic, unitary, supercuspidal representation of $\G(n,F)$ with central character 
$\omega = \omega_\pi.$  Let 
$I(s) = \operatorname{Ind}\left(|\det|^s \sigma \otimes \pi \right)$ denote the parabolically 
induced representation of $\G(m+n,F).$  

If $n=0$ and $m$ is odd, i.e., non-self-associate parabolic (cf. Remark \ref{non-self-assoc}), 
then $I(s)$ is always irreducible.  Otherwise, the following hold: 

\begin{itemize}

\item If $\sigma$ is not $\omega^{-1}$-self-dual, then $I(s)$ is irreducible for $s \ge 0.$ 

\item If $\sigma$ is $\omega^{-1}$-self-dual and $L(s, \sigma, \wedge^2\otimes \omega^{-1})$ 
has a pole at $s=0,$ then $I(s)$ 
is irreducible for $ 0 \le s < 1/2,$ reducible for $s=1/2$, and irreducible for $s > 1/2.$ 
(The reducibility point is the same whether $n = 0$ or $n \ge 1.$) 

\item If $\sigma$ is $\omega^{-1}$-self-dual, $L(s, \sigma, \sym^2 \otimes \omega^{-1})$ 
has a pole at $s=0$ and $n=0,$ 
then $I(s)$ is reducible for $s=0$ and irreducible for $s > 0.$ 

\item If $\sigma$ is $\omega^{-1}$-self-dual, $L(s, \sigma, \sym^2\otimes \omega^{-1})$ 
has a pole at $s=0,$ 
$n \ge 1,$ and 
$\sigma$ appears as a component of the transfer of $\pi$ to $\GL(2n,F),$  then $I(s)$ 
is irreducible for $0 \le s < 1,$ reducible for $s=1,$ and irreducible for $s >1.$ 

\item If $\sigma$ is $\omega^{-1}$-self-dual, $L(s, \sigma, \sym^2\otimes \omega^{-1})$ 
has a pole at $s=0,$ $n \ge 1$ and 
$\sigma$ does not appear as a component of the transfer of $\pi,$ then $I(s)$ 
is reducible for $s=0$ and irreducible for $s > 0.$ 

\end{itemize}

\end{prop}

\begin{prop}\label{unitary} Reducibility for $\G(n) = \Uef(2n)$ and $\Uef(2n+1).$ \newline
Let $m \ge 1$ and let $\sigma$ be an irreducible, unitary, supercuspidal 
representation of $\GL(m, E).$  Let $n \ge 0$ and let $\pi$ be an irreducible, 
generic, unitary, supercuspidal representation of $\G(n,F).$  We may consider 
$\sigma \otimes \pi$ as a representation of $\M(F),$ where 
$\M \cong \operatorname{Res}_{E/F} \GL(m) \times \G(n)$ is a maximal 
Levi subgroup in $\G(m+n).$ Let 
$I(s) = \operatorname{Ind}\left(|\det|^s \sigma \otimes \pi \right)$ denote the parabolically 
induced representation of $\G(m+n,F).$  The following hold: 

\begin{itemize}

\item If $\sigma$ is not conjugate-self-dual, then $I(s)$ is irreducible for $s \ge 0.$ 

\item If $\sigma$ is conjugate-self-dual and 
$L(s, \sigma, r_A)$ has a pole at $s=0$ when $\G(n) = \Uef(2n)$ or  
$L(s, \sigma, r_A \otimes \delta_{E/F})$ has a pole at $s=0$ when $\G(n) = \Uef(2n+1),$
then $I(s)$ is irreducible for $ 0 \le s < 1/2,$ reducible for $s=1/2$, and irreducible for $s > 1/2.$ 
(The reducibility point is the same whether $n = 0$ or $n \ge 1.$) 

\item If $\sigma$ is conjugate-self-dual,
$L(s, \sigma, r_A \otimes \delta_{E/F})$ has a pole at $s=0$ when $\G(n) = \Uef(2n)$ or  
$L(s, \sigma, r_A)$ has a pole at $s=0$ when $\G(n) = \Uef(2n+1),$
and $n=0,$ 
then $I(s)$ is reducible for $s=0$ and irreducible for $s > 0.$ 

\item If $\sigma$ is conjugate-self-dual, 
$L(s, \sigma, r_A \otimes \delta_{E/F})$ has a pole at $s=0$ when $\G(n) = \Uef(2n)$ or  
$L(s, \sigma, r_A)$ has a pole at $s=0$ when $\G(n) = \Uef(2n+1),$
$n \ge 1,$ 
and 
$\sigma$ appears as a component of the transfer of $\pi$ 
to $\GL(2n,F)$ when $\G(n) = \Uef(2n)$ or  
to $\GL(2n+1,F)$ when $\G(n) = \Uef(2n+1),$ 
then $I(s)$ 
is irreducible for $0 \le s < 1,$ reducible for $s=1,$ and irreducible for $s >1.$ 

\item If $\sigma$ is conjugate-self-dual, 
$L(s, \sigma, r_A \otimes \delta_{E/F})$ has a pole at $s=0$ when $\G(n) = \Uef(2n)$ or  
$L(s, \sigma, r_A)$ has a pole at $s=0$ when $\G(n) = \Uef(2n+1),$
$n \ge 1,$ 
and 
$\sigma$ does not appear as a component of the transfer of $\pi,$ then $I(s)$ 
is reducible for $s=0$ and irreducible for $s > 0.$ 

\end{itemize}

\end{prop}

We remark that some of the $L$-function conditions in the above propositions can also be 
restated in other ways such as conditions on the parity of the integer $m.$

\end{document}